\documentclass[ a4paper, oneside]{amsart2000}

\RequirePackage{mathmacros}
\RequirePackage{enumerate}
\usepackage{xr-hyper}
\usepackage[bookmarks=true,  colorlinks=true, pagebackref]{hyperref}

\newcommand{\anl}{\htmladdnormallink}

\theoremstyle{plain}

\theoremstyle{definition}

\swapnumbers
\theoremstyle{remark}

\numberwithin{equation}{section}

\begin{document}
\title{Estimates for the volume of a Lorentzian manifold}

\author{Claus Gerhardt}
\address{Ruprecht-Karls-Universit\"at, Institut f\"ur Angewandte Mathematik,
Im Neuenheimer Feld 294, 69120 Heidelberg, Germany}
\email{gerhardt@math.uni-heidelberg.de}
\urladdr{\url{http://www.math.uni-heidelberg.de/studinfo/gerhardt/}}

%
\date{April 18, 2002}
%


\begin{abstract} We prove new estimates for the volume of a Lorentzian
manifold and show especially that cosmological spacetimes with crushing
singularities have finite volume.
\end{abstract}
\maketitle
\thispagestyle{empty}

\setcounter{section}{-1}
\section{Introduction} 

\cvb
Let $N$ be a  $(n+1)$-dimensional Lorentzian manifold and suppose that $N$ can be
decomposed in the form

\begin{equation}\lae{0.1}
N=N_0\uu N_-\uu N_+,
\end{equation}

\cvm 
\nd where $N_0$ has finite volume and $N_-$ resp. $N_+$ represent the critical
past resp. future Cauchy developments with not necessarily a priori bounded
volume. We assume that $N_+$ is the future Cauchy development of a Cauchy
hypersurface $M_1$, and $N_-$ the past Cauchy development of a hypersurface
$M_2$, or, more precisely, we assume the existence of a time function $x^0$,
such that

\begin{equation}
\begin{aligned}
N_+&={x^0}^{-1}([t_1,T_+)),&\qq M_1=\{x^0=t_1\}&,\\
N_-&={x^0}^{-1}((T_-,t_2]),&\qq M_2=\{x^0=t_2\}&,
\end{aligned}
\end{equation}

\cvm
\nd and that the Lorentz metric can be expressed as

\begin{equation}\lae{0.3}
d\bar s^2=e^{2\psi}\{-{dx^0}^2+\s_{ij}(x^0,x)dx^idx^j\},
\end{equation}

\cvm
\nd where $x=(x^i)$ are local coordinates for the space-like hypersurface $M_1$
if $N_+$ is considered resp. $M_2$ in case of $N_-$.

The coordinate system $(x^\al)_{0\le\al\le n}$ is supposed to be future
directed, i.e. the \tit{past} directed unit normal $(\nu^\al)$ of the level sets

\begin{equation}
M(t)=\{x^0=t\}
\end{equation}

\cvm
\nd is of the form

\begin{equation}\lae{0.5}
(\nu^\al)=-e^{-\psi}(1,0,\ldots,0).
\end{equation}

\cvm
If we assume the mean curvature of the slices $M(t)$ with respect to the past
directed normal---cf. \ci[Section 2]{cg8} for a more detailed explanation of our
conventions---is strictly bounded away from zero, then, the following volume
estimates can be proved

\bt\lat{0.1}
Suppose there exists a positive constant $\e_0$ such that

\begin{align}
H(t)&\ge \e_0&\A\,t_1\le t< T_+&,\lae{0.6}\\
\intertext{and}
H(t)&\le-\e_0&\A\,T_-<t\le t_2&,\lae{0.7}
\end{align}

\cvm
\nd then

\begin{align}
\abs{N_+}&\le \frac1{\e_0}\abs{M(t_1)},\\
\intertext{and}
\abs{N_-}&\le \frac1{\e_0}\abs{M(t_2}.
\end{align}

These estimates also hold locally, i.e. if $E_i\su M(t_i)$, $i=1,2$, are measurable
subsets and $E_1^+,E_2^-$ the corresponding future resp. past directed
cylinders, then,

\begin{align}
\abs{E_1^+}&\le\frac1{\e_0}\abs{E_1},\lae{0.10}\\
\intertext{and}
\abs{E_2^-}&\le\frac1{\e_0}\abs{E_2}.
\end{align}
\et

\cvb
\section{Proof of \rt{0.1}}\las{1}

\cvb
In the following we shall only prove the estimate for $N_+$, since the other case
$N_-$ can easily be considered as a future development by reversing the time
direction.

\cvm
Let $x=x(\xi)$ be an embedding of a space-like hypersurface and $(\nu^\al)$ be
the past directed normal. Then, we have the Gau{\ss} formula

\begin{equation}
x^\al_{ij}=h_{ij}\nu^\al.
\end{equation}

\cvm
\nd where $(h_{ij})$ is the second fundamental form, and the Weingarten equation

\begin{equation}
\nu^\al_i=h^k_ix^\al_k.
\end{equation}

\cvm
We emphasize that covariant derivatives, indicated simply by indices, are
always \tit{full} tensors.

\cvm
The slices $M(t)$ can be viewed as  special embeddings of the form

\begin{equation}
x(t)=(t,x^i),
\end{equation}

\cvm
\nd where $(x^i)$ are coordinates of the \tit{initial} slice $M(t_1)$. Hence, the
slices $M(t)$ can be considered as the solution of the evolution problem

\begin{equation}\lae{1.4}
\dot x=-e^\psi \nu, \qq t_1\le t<T_+,
\end{equation}

\cvm
\nd with initial hypersurface $M(t_1)$, in view of \re{0.5}.

\cvm From the equation \re{1.4} we can immediately derive evolution equations
for the geometric quantities $g_{ij}, h_{ij}, \nu$, and $H=g^{ij}h_{ij}$ of $M(t)$, cf.
e.g.
\ci[Section 4]{cg4}, where the corresponding evolution equations are derived in
Riemannian space.

\cvm
For our purpose, we are only interested in the evolution equation for the metric,
and we deduce

\begin{equation}
\dot g_{ij}=\spd{\dot x_i}{x_j}+\spd{x_i}{\dot x_j}=- 2e^\psi h_{ij},
\end{equation}

\cvm
\nd in view of the Weingarten equation.

\cvm
Let $g=\det(g_{ij})$, then,

\begin{equation}\lae{1.6}
\dot g= g g^{ij}\dot g_{ij}=-2e^\psi H g,
\end{equation}

\cvm
\nd and thus, the volume of $M(t), \abs{M(t)}$, evolves according to

\begin{equation}\lae{1.7}
\frac d{dt}  \abs{M(t)}=\int_{M(t_1)}\frac d{dt}\sqrt g=-\int_{M(t)}e^\psi H,
\end{equation}

\cvm
\nd where we shall assume without loss of generality that $\abs{M(t_1}$ is finite,
otherwise, we replace $M(t_1)$ by an arbitrary measurable subset of $M(t_1)$
with finite volume.

\cvm
Now, let $T\in [t_1, T_+)$ be arbitrary and denote by $Q(t_1,T)$ the
cylinder

\begin{equation}\lae{1.8}
Q(t_1,T)=\set{(x^0,x)}{t_1\le x^0\le T},
\end{equation}

\cvm
\nd then,

\begin{equation}\lae{1.9}
\abs{Q(t_1,T)}=\int_{t_1}^T\int_Me^\psi,
\end{equation}

\cvm
\nd where we omit the volume elements, and where, $M=M(x^0)$.

\cvm
By assumption, the mean curvature $H$ of the slices is bounded from below by
$\e_0$, and we conclude further, with the help of \re{1.7},

\begin{equation}
\begin{aligned}
\abs{Q(t_1,T)}&\le\frac 1{\e_0} \int_{t_1}^T\int_Me^\psi H\\
&=\frac1{\e_0}\{\abs{M(t_1)}-\abs{M(T)}\}\\
&\le \frac1{\e_0}\abs{M(t_1)}.
\end{aligned}
\end{equation}

\cvm
Letting $T$ tend to $T_+$ gives the estimate for $\abs {N_+}$.

\cvm
To prove the estimate \re{0.10}, we simply replace $M(t_1)$ by $E_1$.

\cvb
If we relax the conditions \re{0.6} and \re{0.7} to include the case $\e_0=0$, a
volume estimate is still possible.

\cvm
\bt
If the assumptions of \rt{0.1} are valid with $\e_0=0$, and if in addition the
length of any future directed curve starting from $M(t_1)$ is bounded by a
constant $\ga_1$ and the length of past any directed curve starting from $M(t_2)$
is bounded by a constant $\ga_2$, then,
\begin{align}
\abs{N_+}&\le \ga_1\abs{M(t_1)}\\
\intertext{and}
\abs{N_-}&\le \ga_2\abs{M(t_2)}.
\end{align}
\et

\cvm
\bp
As before, we only consider the estimate for $N_+$.

\cvm
From \re{1.6} we infer that the volume element of the slices $M(t)$ is decreasing
in $t$, and hence,
\begin{equation}\lae{1.13}
\sqrt{g(t)}\le \sqrt{g(t_1)}\qq\A\,t_1\le t.
\end{equation}

\cvm
Furthermore, for fixed $x\in M(t_1)$ and $t>t_1$
\begin{equation}\lae{1.14}
\int_{t_1}^te^\psi\le \ga_1
\end{equation}
because the left-hand side is the length of the future directed curve
\begin{equation}
\ga(\tau)=(\tau,x)\qq t_1\le\tau\le t.
\end{equation}

\cvm
Let us now look at the cylinder $Q(t_1,T)$ as in \re{1.8} and \re{1.9}. We have
\begin{equation}
\begin{aligned}
\abs{Q(t_1,T)}&=\int_{t_1}^T\int_{M(t_1)}e^\psi\sqrt{g(t,x)}\le
\int_{t_1}^T\int_{M(t_1)}e^\psi\sqrt{g(t_1,x)}\\[\cma]
&\le \ga_1\int_{M(t_1)}\sqrt{g(t_1,x)}=\ga_1\abs{M(t_1)}
\end{aligned}
\end{equation}
by applying Fubini's theorem and the estimates \re{1.13} and \re{1.14}.
\ep

\cvb
\section{Cosmological spacetimes}\las{2}

\cvb
A cosmological spacetime is a globally hyperbolic Lorentzian manifold $N$ with
compact Cauchy hypersurface $\so$, that satisfies the timelike convergence
condition, i.e.

\begin{equation}
\bar R_{\al\bet}\nu^\al\nu^\bet\ge 0 \qq \A\,\spd\nu\nu=-1.
\end{equation}

\cvm
If there exist crushing singularities, see \ci{es} or \ci{cg1} for a definition, then,
we proved in
\ci{cg1} that
$N$ can be foliated by spacelike hypersurfaces $M(\tau)$ of constant mean
curvature $\tau$, $-\un<\tau<\un$,

\begin{equation}
N=\uuu_{0\ne\tau\in \R[]}M(\tau)\uu{\msc C}_0,
\end{equation}

\cvm
\nd where $\msc C_0$ consists either of a single maximal slice or of a whole
continuum of maximal slices in which case the metric is stationary in $\msc
C_0$. But in any case $\msc C_0$ is a compact subset of $N$.

\cvm
In the complement of $\msc C_0$ the mean curvature function $\tau$ is a regular
function with non-vanishing gradient that can be used as a new time function, cf.
\ci{cg6} for a simple proof.

\cvm
Thus, the Lorentz metric can be expressed in Gaussian coordinates $(x^\al)$ with
$x^0=\tau$ as in \re{0.3}. We choose arbitrary $\tau_2<0<\tau_1$ and de\-fine

\begin{equation}
\begin{aligned}
N_0&=\set{(\tau,x)}{\tau_2\le\tau \le \tau_1},\\
N_-&=\set{(\tau,x)}{-\un<\tau \le \tau_2},\\
N_+&=\set{(\tau,x)}{\tau_1\le \tau<\un}.
\end{aligned}
\end{equation}

\cvm
Then, $N_0$ is compact, and the volumes of $N_-, N_+$ can be estimated by

\begin{align}
\abs{N_+}&\le \frac1{\tau_1}\abs{M(\tau_1)},\\
\intertext{and}
\abs{N_-}&\le \frac1{\abs{\tau_2}}\abs{M(\tau_2)}.
\end{align}

\cvm
Hence, we have proved

\bt
A cosmological spacetime $N$ with crushing singularities has finite volume.
\et

\cvb
\br
Let $N$ be a spacetime with compact Cauchy hypersurface and suppose that a
subset
$N_-\su N$ is foliated by constant mean curvature slices $M(\tau)$ such that

\begin{equation}
N_-=\uuu_{0<\tau\le \tau_2}M(\tau)
\end{equation}

\cvm
\nd and suppose furthermore, that $x^0=\tau$ is a time function---which will be
the case if the timelike convergence condition is satisfied---so that the metric
can be represented in Gaussian coordinates $(x^\al)$ with $x^0=\tau$.

\cvm
Consider the cylinder $Q(\tau,\tau_2)=\{\tau\le x^0\le \tau_2\}$ for some
fixed $\tau$. Then, 

\begin{equation}
\abs{Q(\tau,\tau_2)}=\int_\tau^{\tau_2}\int_Me^\psi=\int_\tau
^{\tau_2}H^{-1}\int_MH e^\psi,
\end{equation}

\cvm
\nd and we obtain in view of \re{1.7}

\begin{equation}
\tau^{-1}_2\{\abs {M(\tau)}-\abs{M(\tau_2)}\}\le\abs{Q(\tau,\tau_2)},
\end{equation}

\cvm
\nd and conclude further

\begin{equation}
\lim_{\tau\ra 0}\msp[2]\abs{M(\tau)}\le \tau_2\abs{N_-}+\abs{M(\tau_2)},
\end{equation}

\nd i.e.

\begin{equation}
\lim_{\tau\ra 0}\msp[2]\abs{M(\tau)}=\un\im \abs{N_-}=\un.
\end{equation}
\er

\cvb
\section{The Riemannian case}

\cvb
Suppose that $N$ is a Riemannian manifold that is decomposed as in \re{0.1} with
metric

\begin{equation}
d\bar s^2=e^{2\psi}\{{dx^0}^2+\s_{ij}(x^0,x)dx^idx^j\}.
\end{equation}

\cvm
The Gau{\ss} formula and the Weingarten equation for a hypersurface now have
the form

\begin{align}
x^\al_{ij}&=-h_{ij}\nu^\al,\\
\intertext{and}
\nu^\al_i&=h^k_ix^\al_k.
\end{align}

\cvm
As default normal vector---if such a choice is possible---we choose the outward
normal, which, in case of the coordinate slices $M(t)=\{x^0=t\}$ is given by

\begin{equation}
(\nu^\al)=e^{-\psi}(1,0,\ldots,0).
\end{equation}

\cvm
Thus, the coordinate slices are solutions of the evolution problem

\begin{equation}
\dot x=e^\psi \nu,
\end{equation}

\cvm
\nd and, therefore,

\begin{equation}
\dot g_{ij}=2e^\psi h_{ij},
\end{equation}

\cvm
\nd i.e. we have the opposite sign compared to the Lorentzian case leading to

\begin{equation}
\frac d{dt}\abs{M(t)}=\int_Me^\psi H.
\end{equation}

\cvm
The arguments in \rs{1} now yield

\bt
\tup{(i)} Suppose there exists a positive constant $\e_0$ such that the mean
curvature $H(t)$ of the slices $M(t)$ is estimated by

\begin{align}
H(t)&\ge \e_0&\A\,t_1\le t< T_+&,\\
\intertext{and}
H(t)&\le-\e_0&\A\,T_-<t\le t_2&,
\end{align}

\cvm
\nd then

\begin{align}
\abs{N_+}&\le \frac1{\e_0}\lim_{t\ra T_+}\abs{M(t)},\\
\intertext{and}
\abs{N_-}&\le \frac1{\e_0}\lim_{t\ra T_-}\abs{M(t}.
\end{align}

\cvm
\tup{(ii)} On the other hand, if the mean curvature $H$ is negative in $N_+$ and
positive in $N_-$, then, we obtain the same estimates as \rt{0.1}, namely,

\begin{align}
\abs{N_+}&\le \frac1{\e_0}\abs{M(t_1)},\\
\intertext{and}
\abs{N_-}&\le \frac1{\e_0}\abs{M(t_2)}.
\end{align}
\et

\cvb

\end{document}